\documentclass[12pt,a4paper,twoside,leqno]{article}
\usepackage[english]{babel}
\usepackage{amssymb}
\usepackage{inputenc}
\newcommand{\ba}{\begin{array}}
\newcommand{\ea}{\end{array}}
\usepackage{amssymb}
\newtheorem{theorem}{Theorem}[section]

\newtheorem{definition}[theorem]{Definition}
\newtheorem{proposition}[theorem]{Proposition}

\newtheorem{corollary}[theorem]{Corollary}

\usepackage{cite}
\DeclareSymbolFont{rsfs}{U}{rsfs}{m}{n}
\DeclareSymbolFontAlphabet{\rsfsmathscr}{rsfs}
\usepackage[mathscr]{eucal}
\usepackage{amscd}
\usepackage{floatflt}
\usepackage{amsmath}
\usepackage{amssymb}

\def\e{{\bf 1}\!\!{\rm I}}
\numberwithin{equation}{section}
\begin{document}

\author{S.Albeverio $^{1},$ Sh.A.Ayupov $^2$, A.A.Rakhimov, $^{2,  *}$, R.A.Dadakhodjaev $^2$}

\title{\bf  Index Theory for Real Factors.}

\date{}

\maketitle

\begin{abstract}
The notion of index for arbitrary real factors is introduced and investigated.
The main tool in our approach is the reduction of real factors to involutive *-anti-automorphisms of their complex enveloping von Neumann algebras.
Similar to the complex case the values of the index for real factors are calculated.
\end{abstract}
\medskip
\medskip

$^1$ Institut f\"{u}r Angewandte Mathematik, Universit\"{a}t Bonn,
Endenicher Allee 60, D-53115 Bonn (Germany); SFB 611; HCM;  BiBoS; IZKS; CERFIM
(Locarno); Acc. Arch. (USI),

e-mail address:\emph{albeverio@uni-bonn.de}

\medskip
$^2$ Institute of Mathematics and information  technologies,
Uzbekistan Academy of Sciences, Dormon Yoli str. 29, 100125,
Tashkent (Uzbekistan),

The Abdus Salam International Centre for Theoretical Physics, Trieste, Italy

e-mail address: \emph{sh\_ayupov@mail.ru}
\medskip



\medskip



\medskip

* Corresponding author. \emph{rakhimov@ktu.edu.tr}

\section{{\large Introduction.} }

In \cite{K1}, H.Kosaki extended the notion of index to an arbitrary (normal faithful) expectation from a factor onto a
subfactor. While Jones' definition of the index is based on the coupling constant, Kosaki's definition of the index of an
expectation relies on the notion of spatial derivatives due to A.Connes \cite{C} as well as the theory of operator-value weights
due to U.Haagerup \cite{H}. In \cite{K1,K2}, it was shown that many fundamental properties of the Jones index in the type II$_1$ case
can be extended to the general setting.
At the present time, the theory of index thanks to works by V.Jones, P.Loi, R.Longo, H.Kosaki and
other mathematicians is well developed and has many applications
in the theory of operator algebras and physics (see also \cite{PL,RL}).\\[-2mm]

In parallel with the theory of an index of complex subfactors the theory of an index of real subfactors has also been intensively developed.
In the papers \cite{AAR1, AAR2, AR1, RK1, R2, RK2} real analogues of Jones' theory of index are considered. In particular,
the notions of the real coupling constant and the index for finite real factors was introduced and investigated.
In the present paper the notion of real index for arbitrary real factors is introduced and investigated.
The main tool in our approach is the reduction of real factors to involutive *-anti-automorphisms of their complex enveloping von Neumann algebras.
Similar to the complex case the values of the index for real factors are calculated.

\section{{\large Preliminaries.} }

Let $B(H)$ be the algebra of all bounded linear operators on a complex Hilbert space $H$. A weakly closed *-subalgebra
$\mathfrak{A}$ containing the identity operator $\e$ in $B(H)$ is called a W$^*$-{\it algebra}. A real *-subalgebra $R\subset B(H)$ is called
a {\it real} W$^*$-{\it algebra} if it is closed in the weak operator topology and \ $R\cap iR=\{0\}$. A real
W$^*$-algebra $R$ is called a {\it real factor} if its center $Z(R)$ consists of the elements $\{\lambda\e , \lambda\in\mathbb{R}\}$.
We say that a real W$^*$-algebra $R$ is of the type I$_{fin}$, I$_{\infty}$, II$_1,$ II$_{\infty},$ \ or \
III$_{\lambda}$, $(0 \leq \lambda \leq 1)$ \ if the enveloping W$^*$-algebra $\mathfrak{A}(R)$ has the corresponding type in the
ordinary classification of W$^*$-algebras. A linear mapping $\alpha$ of an algebra into itself with $\alpha(x^*)=\alpha(x)^*$ \ is called
an *-{\it automorphism} if $\alpha(xy)=\alpha(x)\alpha(y)$; it is called an {\it involutive *-antiautomorphism} if $\alpha(xy)=\alpha(y)\alpha(x)$
and $\alpha ^2(x)=x$. If $\alpha$ is an involutive *-antiautomorphism of a W$^*$-algebra $M$, we denote by $(M,\alpha)$
the real W$^*$-algebra generated by $\alpha$, i.e. $(M,\alpha)=\{x\in M: \ \alpha(x)=x^*\}$. Conversely, every real W*-algebra $R$ is of
the form $(M,\alpha)$, where $M$ is the complex envelope of $R$ and $\alpha$ is an involutive *-antiautomorphism of $M$ (see \cite{ARU,G,R}).
Therefore we shall identify from now on the real von Neumann algebra $R$ with the pair $(M,\alpha)$.

\section{\large Extended positive part of a real W*-algebra.}

We recall the definitions of the extended positive part ${\hat N}^+$ of a W*-algebra $N$ and an operator-valued trace
on a W*-algebra $M^+$ with range in the extended positive part ${\hat N}^+$ of a W*-subalgebra $N\subset M$ \cite{H}.
Let $N^+_*$ be the set of all normal positive linear functionals on $N$. We consider the set ${\hat N}^+$ of all
positively homogeneous additive lower semicontinuous functions $m: N^+_*\to [0,+\infty]$ and embed the cone $N^+$ in
${\hat N}^+$ by identifying an element $x\in N^+$ with the function $m_x$ defined by the relation $m_x(f)=f(x)$,
for all $f\in N^+_*$. For an unbounded self-adjoint positive operator $x$ affiliated with $N$ we denote its
support by $e$ and define the corresponding function $m_x$ by the formula
$$
m_x(f) = \sum f({\overline e}_n x) + (+\infty)f(\e - e) ,
$$
where $f$ is an arbitrary functional in $N^+_*$ and ${\overline e}_n = e_{[n-1,n]}$ are the spectral projections of $x$
corresponding to the set $[n-1,n]$, $n\in \{1,2,\ldots ,\infty\}$. It was shown in \cite{H} that for each $m\in {\hat N}^+$
there exists a positive self-adjoint (but not necessarily bounded) operator $A$ affiliated with $N$ such that
$m=m_A$.

Let $M$ be a W*-algebra, and let $N$ be a W*-subalgebra of $M$. By an {\it operator-valued weight} on the W*-algebra
$M$ with range in ${\hat N}$, or an $N$-{\it valued weight}, we mean a linear map $E:M\to {\hat N}$ such that
$E(yxy^*)=yE(x)y^*$, for $x\in M^+$ and $y\in N$ \cite{H}. The properties of being normal, faithful, or
semifinite are defined for $E$ similarly to the case of linear functionals.

Let $M_E=\{x\in M: \|E(x^*x)\|<\infty\}=\{x\in M: E(x^*x)\in N^+\}$. As is known, $M_E$ is a facial subalgebra of $M$
and an $N$-bimodule (that is, $N\cdot M_E\cdot N \subset M_E$ and $E(yxz)=yE(x)z$, for all $x\in M_E$ and $y,z\in N$);
moreover, $E$ can be uniquely extended to a linear map $E:M_E\to N$, and $E(M_E)$ is a two-sided ideal of $N$.
Hence, $E$ can be uniquely extended to the map ${\overline E}: {\hat M}^+ \to {\hat N}^+$ (see \cite{H}).\\[-2mm]

Now, we recall the definitions of the extended positive part ${\hat R}^+$ of a real W*-algebra $R=(M,\alpha)$ \cite{USh1}.
Let $(R^+_*)_r$ be the set of all normal positive linear functionals on $R$ that vanish on skew-symmetric elements of $R$;
each of these functionals has a unique $\alpha$-invariant normal extension to $M$ (see \cite{A1}).
On the set $(R^+_*)_r$ we consider the family ${\hat R}^+$ of all positively homogeneous additive lower semicontinuous
functions $m: (R^+_*)_r\to [0,+\infty]$; we can embed the cone $R^+$ in ${\hat R}^+$ by identifying an element $x\in R^+$
with the function $m_x$ such that $m_x(f)=f(x)$, for all $f\in (R^+_*)_r$. On the other hand, if $x$ is an unbounded
self-adjoint positive operator affiliated with $R$ and with support $e$, then we define the corresponding function $m_x$ by the formula
$$
m_x(f) = \sum f({\overline e}_n x) + (+\infty)f(\e - e) ,
$$
where $f$ is an arbitrary functional in $(R^+_*)_r$ and ${\overline e}_n = e_{[n-1,n]}$ are the spectral projections of $x$
corresponding to the set $[n-1,n]$, $n\in \{1,2,\ldots ,\infty\}$.

\begin{definition} \label{proposition11}
The extended positive part of a real W*-algebra $R$ is the set ${\hat R}^+$.
\end{definition}

\begin{proposition}{\rm \cite[Proposition 3.1.]{USh1}} For each $m\in {\hat R}^+$ there exist a projection $e\in K$ and
a positive self-adjoint (not necessarily bounded) operator $A$ on $eH$ affiliated with $R$ such that $m=m_A$.
\end{proposition}

It follows from Proposition \ref{proposition11} that ${\hat R}^+ \subset {\hat M}^+$, where $M=\mathfrak{A}(R)$ is the enveloping W*-algebra of $R$.\\[-2mm]

By an {\it operator-valued weight} on a real W*-algebra $R$ with range in the extended positive part ${\hat Q}$ of a real W*-subalgebra $Q$ of $R$
(or a {\it $Q$-valued weight}) we mean a linear map $T:R^+\to {\hat Q}^+$ such that $T(yxy^*)=yT(x)y^*$, for $x\in R^+$ and $y\in Q$ \cite{USh1}.
The properties of being normal, faithful, or semifinite are defined for $T$ similarly to the case of linear functionals.
For completeness let us recall the definitions: an operator-valued weight $T$ on $R$ is said to be\\[-2mm]

- {\it normal} if the convergence $x_i\nearrow x$ implies $T(x_i)\nearrow T(x)$, \ for $x_i,x\in R^+$;\\[-2mm]

- {\it faithful} if the equality $T(x^*x)=0$ implies $x=0$;\\[-2mm]

- {\it semifinite} if the set $R_T=\{x\in R: \|T(x^*x)\|<\infty\}$ is ultraweakly dense in $R$.\\[-2mm]

\noindent The set $R_T$ is also a facial subalgebra of $R$ and $R_T$ is a real $Q$-bimodule;
moreover, $T$ can be uniquely extended to a linear map $T:R_T\to Q$, the set $T(R_T)$ is a two-sided ideal of $Q$.
Hence, the map $T$ can be uniquely extended to the map ${\overline T}: {\hat R}^+ \to {\hat Q}^+$.

\section{\large The existence of a normal operator-valued weights.}

Let $M$ be a W*-algebra, and let $N$ be a W*-subalgebra of $M$. The set of normal faithful semi-finite
weights on $M$ is denoted by $P(M)$; the set of normal faithful semi-finite operator-valued weights from $M$ to $N$ is denoted by $P(M,N)$.

As shown by Haagerup, the following results takes place

\begin{theorem}{\rm \cite[Theorem 5.1.]{H}} \label{theorem41}
Let $\psi\in P(M)$ and $\phi\in P(N)$. If $\sigma ^\psi _t(x)=\sigma ^\phi _t(x)$, for any $x\in N$ and $t\in \mathbb{R}$, then there exists a unique
$T\in P(M,N)$, such that $\psi = \phi \circ T$.
\end{theorem}

\begin{theorem}{\rm \cite[Theorem 5.9.]{H}} \label{theorem42}
\ $P(M,N) \not= \emptyset$ \ $\Leftrightarrow$ \ $P(N',M') \not= \emptyset$
\end{theorem}

\noindent Here $M'$ and $N'$ are the commutant of $M$ and $N$, respectively.
Moreover, Haagerup constructed the canonical order-reversing bijection (denoted by $\Phi : E\to E^{-1}$)
between $P(M,N)$ and $P(N',M')$. For a given $E\in P(M,N)$, the canonical bijection $E^{-1}\in P(N',M')$ is characterized
by
$
d(\phi \circ E)/d\psi = d{\phi}/d(\psi \circ E^{-1}).
$
Here $\psi\in P(M')$ and $\phi\in P(N)$, and $E^{-1}$ depends only on $E$.\\[-2mm]

The result similar to theorem \ref{theorem41} is also valid for real W*algebras. Namely the following theorem establishes
necessary and sufficient conditions for the existence of a normal operator-valued weight on $(M,\alpha)$.

\begin{theorem}{\rm \cite[Theorem 4.1.]{USh1}} \label{theorem43} The following conditions are equivalent:
\begin{itemize}
\item[{\rm (1)}] there exists a normal faithful semi-finite $\alpha$-invariant weight $\varphi$ on $M^+$ and
a normal faithful semi-finite $\alpha$-invariant weight $\psi$ on $N^+$ such that
$\sigma ^\varphi _t(x)=\sigma ^\psi _t(x)$, for any $x\in N$ and $t\in \mathbb{R}$;
\item[{\rm (2)}] there exists a unique normal faithful semi-finite operator-valued weight $T$ on $(M,\alpha)$
such that $\varphi = \psi \circ T$.
\end{itemize}
\end{theorem}

\vspace{0.35cm}

It is easy to see that from Theorems \ref{theorem43} and \ref{theorem41} for real W*-algebras $R=(M,\alpha)$ and $Q=(N,\alpha)$ we have following corollary

\begin{corollary} \label{corollary41}
\ $P(R,Q) \not= \emptyset$ \ $\Rightarrow$ \ $P(M,N) \not= \emptyset$.
\end{corollary}

Let us prove the converse implication.

\begin{theorem} \label{theorem44}
$P(M,N) \not= \emptyset$ \ $\Rightarrow$ \ $P(R,Q) \not= \emptyset$
\end{theorem}

\noindent Proof. Let $T_1:M^+\to {\hat N}^+$ be a normal faithful semi-finite operator-valued weight, i.e. $T_1\in P(M,N)$.
We put
$$
T(x) = \frac{1}{2}(T_1(x) + {\overline \alpha} T_1(x)),
$$
where ${\overline \alpha}$ is the extention of $\alpha$ on ${\hat M}^+$.\\[-2mm]

\noindent Since ${\overline \alpha} T = T$, $R^+\subset M^+$ and ${\hat Q}^+ \subset {\hat N}^+$, then for all $x\in R^+$,
$y=T(x)\in {\hat N}^+$ we have ${\overline \alpha}(y) = {\overline \alpha} T(x) = T(x) = y = y^*$, i.e. $y\in {\hat Q}^+$.
Therefore, $T:R^+\to {\hat Q}^+$. The linearity of $T_1$ implies obviously the linearity of $T$. Let $x\in R^+$ and $y\in Q$.
Then according to $T_1(yxy^*)=yT_1(x)y^*$ and ${\overline \alpha}(y^*)=y$ \ we have
$$
T(yxy^*) = \frac{1}{2}(yT_1(x)y^* + {\overline \alpha}(y^*)T_1(x){\overline \alpha}(y)) =
y\Bigl(\frac{1}{2}(T_1(x) + {\overline \alpha} T_1(x))\Bigr)y^* = yT(x)y^*.
$$
Thus, the map $T:R^+\to {\hat Q}^+$ is an operator-valued weight.\\[-2mm]

If $x_i,x\in R^+\subset M^+$ and $x_i \nearrow x$, then by normality of $T_1$ we obtain $T(x_i)\nearrow T(x)$, i.e.
$T$ is also normal. If $T(x^*x)=0$, then $T_1(x^*x)=0$, therefore from the faithfulness of $T_1$ the faithfulness of $T$ follows.
Finally, since $\|T(x^*x)\|<\infty$ $\Leftrightarrow$ $\|T_1(x^*x)\|<\infty$, semifiniteness of $T_1$ implies semifiniteness of $T$.\\[-2mm]

Thus the map $T:R^+\to {\hat Q}^+$, defined as $T=\frac{1}{2}(T_1+{\overline\alpha}T_1)$ is
a normal faithful semi-finite operator-valued weight, i.e. $T\in P(R,Q)$. \ $\Box$\\[-2mm]

Now, by Theorems \ref{theorem41}-\ref{theorem44} and Corollary \ref{corollary41} we obtain the following corollaries

\begin{corollary} \label{corollary42}
\ $P(R,Q) \not= \emptyset$ \ $\Leftrightarrow$ \ $P(M,N) \not= \emptyset$.
\end{corollary}

\begin{corollary} \label{corollary43}
\ $P(R,Q) \not= \emptyset$ \ $\Leftrightarrow$ \ $P(Q',R') \not= \emptyset$.
\end{corollary}
Here $Q'=(N',\alpha ')$, $R'=(M',\alpha ')$ and $\alpha '$ is the involutive *-antiautomorphism of $N'$ defined as
$\alpha '(\cdot )=J \alpha(J \cdot J) J$, where $J:x\to x^*$ is the canonical conjugate linear isometry (see \cite{ARU}).


\section{\large Index for real finite factors.}

Let $F$ ($\subset B(H)$) be a finite (complex or real) factor with the finite commutant $F'$.
The {\it coupling constant} $\dim _F(H)$ of $F$ is defined as tr$_F(E^{F'}_\xi)/$tr$_{F'}(E^F_\xi)$,
where $\xi$ is a non-zero vector in $H$, tr$_A$ denotes the normalized trace and $E^A_\xi$ is
the projection of $H$ onto the closure of the subspace $A\xi$. This definition, in the complex case, is due
to Murray and von Neumann in \cite{MN1}, and in the real case it is introduced in \cite{AAR1, AAR2}.
In both cases it does not depend on $\xi$.\\[-2mm]

It is known \cite[Theorem 6.4.]{AAR1}, that if $M\subset B(H)=B(H_r)+iB(H_r)$ is a finite factor and $(M,\alpha)\subset B(H_r)$,
where $H_r$ is a real Hilbert space with $H_r+iH_r=H$, then
\begin{equation} \label{eq-for2002}
\dim _M(H) = \dim _{(M,\alpha)}(H_r) = \frac{1}{2} \dim _{(M,\alpha)}(H).
\end{equation}
Consider a subfactor $N\subset M$ such that $\alpha(N)\subset N$.
The {\it index} of $N$ in $M$, denoted by $[M:N]$ is defined as $\dim _N(L^2(M))$ (see \cite{J}),
where $L^2(M)$ the completion of $M$ with respect to the norm $\|x\|_2=\tau(x^*x)^{1/2}$.
Similarly, the {\it index} of $(N,\alpha)$ in $(M,\alpha)$, denoted by
$[(M,\alpha):(N,\alpha)]$, or by $[R:Q]$, is defined as $\dim _{(N,\alpha)}(L^2(M,\alpha))$ (see \cite{AAR1, AAR2}).
Between real and complex indices there is the following relation
$$
[(M,\alpha):(N,\alpha)]=[M:N], \ {\rm i.e.} \ [R:Q]=[R+iR:Q+iQ].
$$
(see \cite[Theorem 7.2]{AAR1} and \cite[Theorem 8]{AAR2}).\\[-2mm]

Considering a complex factor $M$ as a real W*-algebra in view of \eqref{eq-for2002} we may put
$[M:(M,\alpha)]=2[(M,\alpha):(M,\alpha)]=2$, i.e. $[M:R]=2$.\\[-2mm]


For example, if $M$ is a factor of type I$_4$, then up to isomorphisms it has seven real W*-subalgebras different from $M$,
which are real or complex subfactors of $M$: \
$\mathbb{R}$, $\mathbb{C}$, $\mathbb{H}$, $M_2(\mathbb{R})$, $M_2(\mathbb{C})$,
$M_2(\mathbb{H})$ and $M_4(\mathbb{R})$, where $\mathbb{H}$ is the quaternion algebra.
The values of the indexes are respectively: \ $[M:M_4(\mathbb{R})]=[M:M_2(\mathbb{H})]=2$, \
$[[M:M_2(\mathbb{C})]=[M_4(\mathbb{R}):M_2(\mathbb{R})]=[M_2(\mathbb{H}):\mathbb{H}]=4$, \
$[M:M_2(\mathbb{R})]=8$, \ $[M:\mathbb{C}]=[M_4(\mathbb{R}):\mathbb{R}]=[M_2(\mathbb{H}):\mathbb{R}]=16$.
\ $[M:\mathbb{R}]=32$.\\[-2mm]

If $M$ is of type I$_5$, then up to isomorphisms it has three real W*-subalgebras different from $M$,
which are real or complex subfactors of $M$: \
$\mathbb{R}$, $\mathbb{C}$, $M_5(\mathbb{R})$. The corresponding indexes are: \ $[M:M_5(\mathbb{R})]=2$,
\ $[M:\mathbb{C}]=[M_5(\mathbb{R}):\mathbb{R}]=25$, \ $[M:\mathbb{R}]=50$.

\vspace{0.5cm}

We have calculated the value of the index in the above examples. It turns out that the index may be calculated also in the general case.
V.Jones in \cite{J} has proved a theorem on the values of the index for subfactors of finite factors.
Let us recall this theorem

\vspace{0.25cm}

\begin{theorem}[\cite{J}, Theorem 4.3.1] \label{teorem5001}
Let $M$ be a finite factor, and let $N$ be a subfactor of $M$ with $[M:N]<\infty$.
Then one has either \ $\displaystyle [M:N]=4\cos ^2\frac{\pi}{q}$ \ for some integer \ $q\geq 3$ \ or \
$[M:N]\geq 4$.
\end{theorem}

From Theorem \ref{teorem5001} we obtain the following real version of the above theorem.

\begin{theorem}[\cite{AAR1}, Theorem 7.5.] \label{teorem5002}
Let $M$ be a finite factor and let $N$ be a subfactor of $M$ with $[M:N]<\infty$.
Given be an involutive *-antiautomorphism $\alpha$ of
$M$ with $\alpha(N)\subset N$, put $R=(M,\alpha)$, $Q=(N,\alpha)$.
Then one has either \ $\displaystyle
[(M,\alpha):(N,\alpha)]=4\cos ^2\frac{\pi}{q}$ \ for some integer \ $q\geq 3$ \ or \ $[(M,\alpha):(N,\alpha)]\geq 4$, i.e. \
$\displaystyle [R:Q]=4\cos ^2\frac{\pi}{q}$ \ for some integer \ $q\geq 3$ \ or \ $[R:Q]\geq 4$.
\end{theorem}

\section{Extension of the index theory to arbitrary real factors.}

Let now $M$ be a $\sigma$-finite factor and let $N$ be a subfactor of $M$ with $\alpha(N)\subset N$.
We recall that a linear positive mapping $E:M\to N$ (or $E:(M,\alpha)\to (N,\alpha)$) is called {\it the conditional expectation} with
respect to the W$^*$-subalgebra N if the following conditions are satisfied:
\begin{itemize}
\item[{\rm (i)}] $E(\e)=\e$;
\item[{\rm (ii)}] $E(E(x)y)=E(x)E(y)=E(xE(y))$;
\item[{\rm (iii)}] $E(x)^*E(x)\leq E(x^*x)$, \ $\forall x,y\in M$.
\end{itemize}
We fix a normal conditional expectation $E$ from $(M,\alpha)$ onto $(N,\alpha)$.
The existence of $E$ follows from \cite[Theorem 1.]{USh2}.
For this it suffices to take a normal faithful semi-finite $\alpha$-invariant weight $\varphi$ on $M$
with $\sigma ^\varphi _t(M)=M$ ($\forall t\in \mathbb{R}$), where
$\sigma ^\varphi _t$ is the modular automorphism group of a weight $\varphi$.
The extension of $E$ on $M$ will be denoted by ${\overline E}$.
Since $E$ is an operator-valued weight by Corollary \ref{corollary43} we get $E^{-1}\in P((N',\alpha '),(M',\alpha '))$.
Similarly by Theorem \ref{theorem42} we have ${\overline E}\in P(M,N)$ and ${\overline E}^{-1}\in P(N',M')$.
By the proof of Theorems \ref{theorem43} (i.e. of \cite[Theorem 4.1.]{USh1}), \ \ref{theorem44} \ and Theorem 1 of \cite{USh2}
we obtain ${\overline E}^{-1}|_{(N',\alpha ')}=E^{-1}$, i.e. ${\overline E}^{-1} = {\overline {E^{-1} }}$.\\[-2mm]

For any unitary $u\in M'$, we have
\begin{equation} \label{eq-for2003}
u{\overline E}^{-1}(\e)u^* = {\overline E}^{-1}(u \e u^*) = {\overline E}^{-1}(\e) .
\end{equation}
It is obvious that ${\overline E}(\e)=E(\e)=\e$, but in general we have ${\overline E}^{-1}(\e)=E^{-1}(\e)\not=\e$.
Since $M$ is a factor, by \eqref{eq-for2003} \ ${\overline E}^{-1}(\e) = E^{-1}(\e)$ is a scalar (possibly $+\infty$).\\[-2mm]

Kosaki defined, in \cite{K1}, the notion of the index as $[M:N]={\overline E}^{-1}(\e)$ and showed that when $M$ is a finite factor, then
his definition coincides with Jones' definition.
Moreover, in Theorem 5.4 \cite{K1} he proved that if $M$ is a $\sigma$-finite factor and $N$ its subfactor, then
similarly to the finite case one has ${\overline E}^{-1}(\e)\in\{4\cos ^2\pi /n (n\geq 3)\}\cup[4,+\infty]$.\\[-2mm]

Now following Kosaki we introduce the following
\begin{definition} \label{definition6001}
The index of $Q=(N,\alpha)$ in $R=(M,\alpha)$, denoted by $[R:Q]$ or by $[(M,\alpha):(N,\alpha)]$,
is defined as the scalar $E^{-1}(\e)$.
\end{definition}
By this definition between real and complex indices there is the following relation
$$
[(M,\alpha):(N,\alpha)] = [M:N], \ {\rm i.e.} \ [R:Q]=[R+iR:Q+iQ].
$$
Thus, we obtain the following real version of the index theorem.

\begin{theorem} \label{teorem6001}
Let $M$ be a $\sigma$-finite factor and let $\alpha$ be an involutive *-antiautomorphism of $M$.
If $N$ is a subfactor of $M$ with $\alpha(N)\subset N$, then
one has either \ $\displaystyle [(M,\alpha):(N,\alpha)]=4\cos ^2\frac{\pi}{q}$ \ for some integer \ $q\geq 3$ \ or \
$[(M,\alpha):(N,\alpha)]\geq 4$.
\end{theorem}

\vspace{2.5cm}

\textbf{Ackowledgments.} \emph{The second, third and the fourth named authors would like to acknowledge the hospitality
of the "Institut f$\ddot{u}$r Angewandte Mathematik",
Universit$\ddot{a}$t Bonn (Germany). This work was supported in part by the DFG AL 214/36-1 project (Germany).}


\newpage

\bigskip


\begin{thebibliography}{99}

\bibitem{AAR1} Albeverio S., Ayupov Sh. A., Rakhimov A. A., Dadakhodjaev R. A. Index for Finite Real Factors.
Preprint of the Universit\"at Bonn (Institut f\"ur Numerische Simulation, Institut f\"ur Angewandte Mathematik,
Mathematisches Institut),  BONN, SFB611, N 458, (2009) 26p. See also: Cornell University, arXiv.org,
Mathematics, Operator Algebras, arXiv: 0902.1152v1.

\bibitem{AAR2} Albeverio S., Ayupov Sh. A., Rakhimov A. A., Dadakhodjaev R. A. On Jones' Index for  Real W*-algebras.
Eurasian Math. J., 1:4 (2010), 5-19.

\bibitem{A1} Ayupov Sh.A. Classification and representation of ordered Jordan algebras. Fan, Tashkent, 1986.

\bibitem{AR1} Ayupov Sh.A., Rakhimov A.A. Real W*-algebras, Actions of groups and Index theory for real factors.
VDM Publishing House Ltd. Beau-Bassin, Mauritius. ISBN 978-3-639-29066-0. 2010, p.138.

\bibitem{ARU} Ayupov Sh.A., Rakhimov A.A. and Usmanov Sh.M., Jordan, Real and Lie Structures in Operator
Algebras, Kluw.Acad.Pub.,MAIA. 418 (1997) 235p.

\bibitem{C} Connes A., Spatial theory of von Neumann algebras, J.Funct. Anal. 35 (1980) 153--164.


\bibitem{G} Giordano T., Antiautomorphismes involutifs des facteurs de von Neumann injectifs,
Th\`{e}se, Universite de Neuch\^{a}tel (1981) 106p.


\bibitem{H} Haagerup U., Operator valued weights in von Neumann algebras I, II,
           J.Funct. Anal. 32 (1979) 175--206; 33 (1979) 339--361.

\bibitem{J} Jones V.F.R., Index for Subfactors, Inventiones Math. 72 (1983) 1--25.


\bibitem{K1} Kosaki H., Extension of Jones' Theory on index to arbitrary factors,
J.Funct. Anal. 66 (1986) 123--140.

\bibitem{K2} Kosaki H., A remark on the minimal index of subfactors,
J.Funct. Anal. 107 (1992) 458--470.



\bibitem{PL} Loi P.H., On the theory of index for type III factors, J. Operator Theory 28 (1992) 251--265.

\bibitem{RL} Longo R., Minimal index and Braided Subfactors, J.Funct. Anal. 109 (1992) 98--112.

\bibitem{MN1} Murray F.J. and von Neumann J., On rings of operators, Ann. Math. 37 (1936) 116--229.



\bibitem{RK1} Rakhimov A.A., Kesicio\u glu Y., Cansu M.N. Index for real factors.
Itogi Nauki. Southern Federal Region. Research on the mathematical analysis. Vol. 3, Vladikavkaz, 2009, pp.200-212.

\bibitem{R2} Rakhimov A.A. Index for finite real factors. Abstracts of the third congress of the world mathematical
society of turkic countries June 30-July 4, Almaty "Kazak universiteti", (2009), Vol. 1, 127p.

\bibitem{RK2} Rakhimov A.A., K\"or T., Kesicio\u glu Y. Index for Finite Real Factors. XX.Ulusal Matematik Sempozyumu,
Bildiri \"ozetleri, 5-8 Eyl\"ul, Atat\"urk \"Universitesi, Erzurum (2007).

\bibitem{R} Rakhimov A.A., Injective real W*-factors of type III$_\lambda$, $0<\lambda <1$,
Funct. Anal. and its Applications 3 (1997) 41--44.



\bibitem{USh1} Usmanov Sh.M. Operator-valued weights on real W*-algebras and reversible JW-algebras.
Sbornik: Mathematics, Volume 190, Number 10, 1999, 105-122.


\bibitem{USh2} Usmanov Sh.M. Conditional expectations on real W*-algebras and JW-algebras.
Izv. Vyssh. Uchebn. Zaved. Mat., 2001, no. 7, 43-47.


\end{thebibliography}
\end{document}